\documentclass[12pt]{article}

\usepackage{amssymb,amsmath}

\newtheorem{theorem}{Theorem}[section]
\newtheorem{corollary}[theorem]{Corollary}
\newtheorem{lemma}[theorem]{Lemma}
\newtheorem{proposition}[theorem]{Proposition}

\title{A direct approach to simplicity of solution manifolds}

\author{Hans-Otto Walther}

\begin{document}
	
\maketitle

\begin{abstract}
\noindent
Differential equations with state-dependent delays define a semiflow of continuously differentiable solution operators in general only on the associated {\it solution manifold} $X\subset C^1([-h,0],\mathbb{R}^n)$. For systems with discrete state-dependent delays we construct a diffeomorphism on a neighbourhood of $X$ which takes $X$ to an open subset of the subspace given by $\phi'(0)=0$. This is in line with earlier work on the nature of solution manifolds. The present approach, however, is new and dismisses all hypotheses beyond smoothness which have been instrumental so far. Compared to a recent case study it is more direct in the sense
that theory of {\it algebraic delay systems} is avoided. 
\end{abstract}

\bigskip

\noindent
Key words: Delay differential equation,  state-dependent delay, solution manifold

\medskip

\noindent
2020 AMS Subject Classification: Primary: 34K43, 34K19, 34K05; Secondary: 58D25.
	
\section{Introduction}

For delay differential equations like 
$$
x'(t)=g(x(t-d))
$$
with a delay $d=d(x(t))\in[0,h]$ which is state-dependent 
linearization is possible on a particular state space, the {\it solution manifold}. When brought to the general form of a retarded functional differential equation
\begin{equation}
x'(t)=f(x_t),
\end{equation} 
where $f:U\to\mathbb{R}^n$ is defined on a set of maps $[-h,0]\to\mathbb{R}^n$ with a left derivative at $t=0$ and $x_t$ denotes the shifted restriction  
$x|[t-h,t]$,
then the solution manifold is the set
$$
X=X_f=\{\phi\in U:\phi'(0)=f(\phi)\}.
$$
Let $C^1_n=C^1([-h,0],\mathbb{R}^n)$ and $C_n=C([-h,0],\mathbb{R}^n)$. For $U\subset C^1_n$ open and $f$ continuously differentiable with the additional property that

\medskip

(e) {\it each derivative $Df(\phi):C^1_n\to\mathbb{R}^n$, $\phi\in U$, has a linear extension $D_ef(\phi):C_n\to\mathbb{R}^n$
and the map
$$
U\times C_n\ni(\phi,\chi)\mapsto D_ef(\phi)\chi\in\mathbb{R}^n
$$
is continuous}

\medskip

the set $X_f\subset C^1_n$ if non-empty  is a continuously differentiable submanifold of codimension $n$, and the segments $x_t$  of
differentiable solutions $x:[-h,t_x)\to\mathbb{R}^n$, $0<t_x\le\infty$, of Eq. (1) define a continuous semiflow on $X_f$ whose solution operators $x_0\mapsto x_t$, $0\le t$, are all continuously differentiable \cite{W1,HKWW}. 

\medskip

The mild smoothness property (e) is a variant of a property introduced in \cite{M-PNP}. From property (e) one finds easily that the map
$$
U\ni\phi\mapsto D_ef(\phi)\in L_c(C_n,\mathbb{R}^n)
$$ 
is locally bounded. However, in general this map is not continuous, which happens in particular for $f$ representing equations with state-dependent delay.

\medskip

In case $f=0$ the solution manifold is the closed subspace
$$
X_0=\{\phi\in C^1_n:\phi'(0)=0\}.
$$
The present paper continues a series of results \cite{KR,W6,KW2} about the nature of solution manifolds in general. They all guarantee, under various hypotheses on $f$, that the solution manifold is simple, or admits {\it regularization}, in the sense that there is a diffeomorphism on a neighbourhood of $X_f$ which transforms $X_f$ to an open subset of the subspace $X_0$. In particular, the solution manifold $X_f$ has a global chart. 

\medskip

Let $k$ delay functionals
$d_{\kappa}:C^1_n\supset U_d\to[-h,0]\subset\mathbb{R}$, $\kappa=1,\ldots,k$, and a map $g:\mathbb{R}^{kn}\to \mathbb{R}^n$ be given. As in \cite{KW2}
we consider the  system of $n$ differential  equations
\begin{equation}
x'(t)=g(x(t+d_1(x_t)),\ldots,x(t+d_k(x_t))) 
\end{equation}
with $k$ state-dependent delays, where
the row vector is an abbreviation for the column vector
$$
(x_1(t+d_1(x_t)),\ldots,x_n(t+d_1(x_t)),\ldots,x_1(t+d_k(x_t)),\ldots,x_n(t+d_k(x_t)))^{tr}\in\mathbb{R}^{kn}.
$$
In order to rewrite the system (2) in the form (1) it is convenient to introduce the functional
$v:U_d\to\mathbb{R}^{kn}$ whose components $v_{\mu}:U_d\to\mathbb{R}$, $\mu=1,\ldots,kn$, are defined by 
$$
v_{\mu}(\phi)=\phi_{\nu}(d_{\kappa}(\phi))
$$ 
with $\kappa\in\{1,\ldots,k\}$ and $\nu\in\{1,\ldots,n\}$ given by
$\mu=(\kappa-1)n+\nu$. We set
$$
U=U_d\cap v^{-1}(V)\subset C^1_n
$$ 
and define $f:U\to\mathbb{R}^n$ by
$$
f(\phi)=g(v(\phi)).
$$
Then Eqs. (2) and (1) coincide for solutions with segments in $U$.

\medskip

The results \cite{KR,W6,KW2} reqire assumptions in addition to smoothness, which are boundedness of $D_ef$ \cite{KR}, or delays being bounded away from zero in the simplest case \cite{W6}, or factorization of the delays which involves a linear map into a finite-dimensional space \cite{W6,KW2}. As an example which escapes such assumptions let us mention 
the protein synthesis model from \cite{GHMWW}.
\medskip

In the present paper we dismiss assumptions beyond smoothness and consider system (2) for

\medskip

(g,d) {\it $V$ and $U_d$ open sets, $g$ continuously differentiable, and all functionals $d_{\kappa}$, $\kappa=1,\ldots,k$, continuously differentiable with property (e).}

\medskip

The main results Theorem 3.4 and Theorem 4.2 guarantee regularization for restrictions $X_f\cap\hat{U}$ of 
the solution manifold to certain open sets $\hat{U}\subset U$.  
In Theorem 3.4, which may be considered the more general result, the sets $\hat{U}$ are convex open subsets of nested sets $U_{\ast}$ which exhaust $U$. Theorem 4.2 asserts regularization 
for every manifold $X_f\cap U_{\ast}$. The restriction by convexity is among others overcome by assuming more smoothness of the delay functionals.
Regularization of the whole solution manifold is easily obtained under
additional boundedness assumptions. 

\medskip

In the final Section 5 below we give examples for which the hypotheses of Theorems 3.4 or 4.2 are satisfied while for each one of the results in \cite{KR,W6,KW2} on simplicity of solution manifolds one or more hypotheses are violated.

\medskip

The present paper has a kind of precursor in the case study \cite{W11} where $n=1=k$. A major difference between both papers is that \cite{W11} uses the framework of algebraic-delay systems \cite{W4,W10} and begins with lifting $U$ to the graph of the delay in $C^1\times\mathbb{R}$, which
causes a considerable detour compared to the present direct approach.  
Another difference is that in \cite{W11} the single delay is assumed to be strictly positive, which is more severe than it may seem at first glance. 

\medskip

An open problem which remains is regularization of the whole solution manifold under the smoothness hypotheses of either Theorem 3.4 or Theorem 4.2.

\medskip

{\bf Notation, conventions.} For subsets $A\subset B$ of a topological space $T$ we say $A$ is open in $B$ if $A$ is open with respect to the relative topology on $B$. The closure and the interior of $A$ are denoted by $cl\,A$ and $int\,A$, respectively. 

\medskip

For elements $v,\tilde{v}$ of a vector space $V$,
$$
v+[0,1](\tilde{v}-v)=\{w\in V: \mbox{There exists}\,\,t\in[0,1]\,\,\mbox{such that}\,\,w=v+t(\tilde{v}-v)\}.
$$

For elements $x$ and nonempty subsets $A$ of normed vector spaces,
$$
dist(x,A)=\inf_{a\in A}|x-a|
$$
Derivatives  and partial derivatives of a map at a given argument are continuous linear maps, indicated by a capital $D$. For maps on domains in $\mathbb{R}^n$, $\partial_{\nu}g(x)=D_{\nu}g(x)1$. In case $n=1$, $\phi'(t)=D\phi(t)1$.

\medskip

Differentiation as a linear operator $C^1_n\to C_n$ is denoted by $\partial$.

\medskip

On $\mathbb{R}^n$ we use a norm with 
$$
\max_{\nu=1,\ldots,n}|x_{\nu}|\le|\xi|\le\sum_{\nu=1}^n|\xi_{\nu}|.
$$ 
For reals $a<b$ we consider the Banach space  $C([a,b],\mathbb{R}^n)$ of continuous maps $[a,b]\to\mathbb{R}^n$ with the norm given by $|\chi|=\max_{a\le t\le b}|\chi(t)|$. Due to the choice of the norm on $\mathbb{R}^n$ we have
$$
\max_{\nu=1,\ldots,n}|\phi_{\nu}|\le|\phi|\le\sum_{\nu=1}^n|\phi_{\nu}|\quad\mbox{on}\quad C_n.
$$
On the space $C^1([a,b],\mathbb{R}^n)$ of continuously differentiable maps $[a,b]\to\mathbb{R}^n$ we use the norm given by $|\phi|_1=|\phi|+|\partial\phi|$. 

\medskip

$I_n$ denotes the continuous linear inclusion map $C^1_n\ni\phi\mapsto\phi\in C_n$, and $I=I_1$.

\medskip

For $\phi\in C_n$ and $\delta>0$, $N_{\delta}(\phi)=\{\chi\in C_n:|\chi-\phi|<\delta\}$.

\medskip

On the space $L_c(B,E)$ of continuous linear maps $L$ from a normed space $B$ into a normed space $E$ we use the norm given by
$$
|L|_{L_c(B,E)}=\sup_{|b|\le1}|Lb|
$$

\medskip

The tangent space of a continuously differentiable submanifold $M$ of a Banach space $B$ at a point $z\in M$ is the vectorspace $T_zM$ of all
tangent vectors $\zeta=Dc(0)1$ of continuously differentiable curves $c:(a,b)\to B$ with $a<0<b$, $c(0)=z$, $c((a,b))\subset M$.

\section{Preliminaries, and a lemma}

{\bf Remarks} (i) {\it If $d:C^1_n\supset U_d\to[-h,0]^k\subset\mathbb{R}^k$ is the restriction of a continuously differentiable functional $\Delta:C_n\supset W_{\Delta}\to[-h,0]^k\subset\mathbb{R}^k$, $W_{\Delta}$ open,  in the sense that}
$$
U_d=I_n^{-1}(W_{\Delta})= W_{\Delta}\cap C^1_n\quad\mbox{and}\quad
d(\phi)=\Delta(I_n\phi)\quad\mbox{on}\quad U_d
$$ 
{\it then $U_d$ is open, $d$ is continuously differentiable and has property (e) with}
$$
D_ed(\phi)\chi=D\Delta(I_n\phi)\chi\quad\mbox{for}\quad\phi\in U_d,\,\,\chi\in C_n,
$$
{\it and the map $U_d\ni\phi\mapsto D_ed(\phi)\in L_c(C_n,\mathbb{R}^n)$ is continuous.}

\medskip

(ii) {\it If all components $f_{\nu}$ of a continuously differentiable map $f:C^1_n\supset U\to\mathbb{R}^n$ have property (e) then also $f$ has property (e), with $D_ef(\phi)$, $\phi\in U$, given by}
$$
(D_ef(\phi))_{\nu}=D_ef_{\nu}(\phi)\quad\mbox{for}\quad\nu=1,\ldots,n,
$$
{\it and conversely.}

\begin{proposition} 
The functional $v$ is continuously differentiable with 
\begin{equation}
Dv_{\mu}(\phi)\hat{\phi}=\hat{\phi}_{\nu}(d_{\kappa}(\phi))+\phi_{\nu}'(d_{\kappa}(\phi))Dd_{\kappa}(\phi)\hat{\phi}
\end{equation} 
for $\mu=(\kappa-1)n+\nu$ with $\kappa\in\{1,\ldots,k\}$ and $\nu\in\{1,\ldots,n\}$, $\phi\in U_d$, and $\hat{\phi}\in C^1_n$. 
Condition (e) is satisfied for $D_ev_{\mu}$ given by Eq. (3) with $\chi\in C_n$ in place of $\hat{\phi}$  and $D_ed_{\kappa}(\phi)\chi$ in place of $Dd_{\kappa}(\phi)\hat{\phi}$.
\end{proposition}

{\bf Proof.} The odd extension of a function $\phi:[-h,0]\to\mathbb{R}$ to a function $\Omega\phi:[-2h,h]\to\mathbb{R}$ by $\Omega\phi(-h-s)=2\phi(-h)-\phi(-h+s)$ and $\Omega\phi(s)=2\phi(0)-\phi(-s)$ for $0<s\le h$
defines a continuous linear map $\Omega:C^1\to C^1([-2h,h],\mathbb{R})$. The evaluation map  $ev:C^1([-2h,h],\mathbb{R})\times(-2h,h)\ni(\psi,s)\mapsto\psi(s)\in\mathbb{R}$ on the open subset $C^1([-2h,h],\mathbb{R})\times(-2h,h)$ of the Banach space $C^1([-2h,h],\mathbb{R})\times\mathbb{R}$  is continuously differentiable, with
$$
D\,ev(\psi,s)(\hat{\psi},\hat{s})=D_1ev(\psi,s)\hat{\psi}+D_2ev(\psi,s)\hat{s}=\hat{\psi}(s)+\psi'(s)\hat{s},
$$
see \cite{W1}. Let $\mu=(\kappa-1)n+\nu$ be given with $\kappa\in\{1,\ldots,k\}$, $\nu\in\{1,\ldots,n\}$. 
Using $v_{\mu}(\phi)=\phi_{\nu}(d_{\kappa}(\phi))=ev(\Omega P_{\nu}\phi,d_{\kappa}(\phi))$, with the projection $P_{\nu}:C^1_n\ni\phi\mapsto\phi_{\nu}\in C^1$, we infer that $v_{\mu}$ is continuously differentiable, with
\begin{eqnarray*}
Dv_{\mu}(\phi)\hat{\phi} & = & D_1ev(\Omega P_{\nu}\phi
,d_{\kappa}(\phi))D(\Omega P_{\nu})(\phi)\hat{\phi}+D_2ev(\Omega P_{\nu}\phi,d_{\kappa}(\phi))Dd_{\kappa}(\phi)\hat{\phi}\\
& = & D_1ev(\Omega P_{\nu}\phi,d_{\kappa}(\phi))\Omega P_{\nu}\hat{\phi}+(\Omega P_{\nu}\phi)'(d_{\kappa}(\phi))Dd_{\kappa}(\phi)\hat{\phi}\\
& = & (\Omega P_{\nu}\hat{\phi})(d_{\kappa}(\phi))+(\Omega P_{\nu}\phi)'(d_{\kappa}(\phi))Dd_{\kappa}(\phi)\hat{\phi}\\  
& = & \hat{\phi}_{\nu}(d_{\kappa}(\phi))+\phi_{\nu}'(d_{\kappa}(\phi))Dd_{\kappa}(\phi)\hat{\phi},
\end{eqnarray*}
which is Eq. (3). In order to verify condition (e) for $v_{\mu}$ use continuity of differentiation $\partial:C^1\to C$, continuity of the evaluation map $ev_C:C\times[-h,0]\ni(\phi,t)\mapsto\phi(t)\in\mathbb{R}$, and property (e) of the delay functional $d_{\kappa}$. $\Box$

\medskip

\begin{proposition} 
The functional $f$ is continuously differentiable and has property (e). If $X_f\neq\emptyset$ then 
$X_f$ is a continuously differentiable submanifold of codimension $n$ in $C^1_n$.
\end{proposition}

\medskip

{\bf Proof.} By the chain rule $f=g\circ (v|U)$ is continuously differentiable, and for every $\nu\in\{1,\ldots,n\}$, $\phi\in U$, $\hat{\phi}\in C^1_n$,
\begin{equation}
Df_{\nu}(\phi)\hat{\phi}=Dg_{\nu}(v(\phi))Dv(\phi)\hat{\phi}=\sum_{\mu=1}^{kn}D_{\mu}g_{\nu}(v(\phi))Dv_{\mu}(\phi)\hat{\phi}
\end{equation} 
Replacing in Eq. (4) the derivatives $Dv_{\mu}(\phi)$ by the extended derivatives $D_ev_{\mu}(\phi):C_n\to\mathbb{R}$ from Proposition 2.1 we obtain  property (e), which in case $X_f\neq\emptyset$ yields the remaining assertion \cite{W1,HKWW}. $\Box$

\medskip

The proofs of Theorems 3.5 and 4.2
below rely on the following lemma which establishes continuously differentiable maps $\tau:{\mathcal V}\to C^1$, ${\mathcal V}\subset\mathbb{R}^m$ open,  with $\tau(y)'(0)=1$ for which all $\tau(y)$ and $D_{\mu}\tau(y)1$ are arbitrarily small in the norm $|\cdot|$ on $C$.

\begin{lemma} 
For every continuous function $H:{\mathcal V}\to(0,\infty)\subset\mathbb{R}$ there exist continuously differentiable maps $\tau:{\mathcal V}\to C^1$ so that for every $y\in {\mathcal V}$,
$$
\tau(y)'(0)=1,\,\,|I\tau(y)|\le H(y),\,\,\mbox{and}\,\, |ID_{\mu}\tau(y)1|\le H(y)\,\,\mbox{for}\,\,\mu=1,\ldots,m.
$$
\end{lemma}


Notice that $\tau(y)$ is transversal to $X_0$ (in case $n=1$).

\medskip

We omit the proof of Lemma 2.3 as it is a simplified version of the proof of \cite[Proposition 2.2]{KW2}, in which
a family of transversals with additional properties (in the kernel of a certain linear map) is constructed.
The function $H$ in Lemma 2.3 corresponds to $h$ in  \cite[Proposition 2.2]{KW2}, and  $\tau$ in Lemma 2.3  corresponds to $H_{\nu}$ in  \cite[Proposition 2.2]{KW2}. The index $\nu$ and \cite[Lemma 2.1]{KW2} prior to \cite[Proposition 2.2]{KW2} are not relevant for Lemma 2.3 and its proof.

\section{Regularization, convex domains}




Maps $U\to C^1_n$ which send a subset $X\subset U$ into $X_0$ can be written as a perturbation of identity
$$
A:U\ni\phi\mapsto\phi-R(\phi)\in C^1_n\quad\mbox{with}\quad R(\phi)'(0)=\phi'(0)\quad\mbox{for}\quad\phi\in X.
$$
The choice $R_{\nu}(\phi)=f_{\nu}(\phi)\tau_{\phi}$ with $\tau_{\phi}\in C^1$ and $(\tau_{\phi})'(0)=1$, like for example 
$$
\tau_{\phi}(t)=t,
$$ 
yields $A(X_f)\subset X_0$, because for every $\phi\in X_f$ and for all $\nu\in\{1,\ldots,n\}$ we have
$$
A_{\nu}(\phi)'(0)=\phi_{\nu}'(0)-f_{\nu}(\phi)\tau_{\phi}'(0)=\phi'(0)-f(\phi)=0.
$$
Also, $A(\phi)=\phi$ on $X_f\cap X_0$. 
All further properties of the map $A$ rely on the choice of the transversals
$\tau_{\phi}\in C^1\setminus X_0,\,\,\phi\in U$.

\medskip

In view of $f(\phi)=g(v(\phi))$ and guided by earlier work \cite{W6,KW2} we consider remainders $R$ whose components have the form
$$
R_{\nu}(\phi)=f_{\nu}(\phi)\tau_{\phi}=g_{\nu}(v(\phi))\tau(v(\phi))\quad\mbox{for}\quad\phi\in U,\,\nu=1,\ldots,n
$$
with a continuously differentiable map $\tau:V\to C^1$ which satisfies $(\tau(y))'(0)=1$ for all $y\in V$, as in Lemma 2.3. The next objective is an analogue of the extension property (e) for such $R$ and an estimate of $D_eR(\phi)$ in terms of $g$ and $d$ and their deivatives.
For $y\in V$ we set
\begin{eqnarray*}
m_g(y) & = & \max_{\nu=1,\ldots,n}|g_{\nu}(y)|,\\
m_{Dg}(y) & = & \max_{\nu=1,\ldots,n;\mu=1,\ldots,kn}|D_{\mu}g_{\nu}(y)1|,
\end{eqnarray*}
and for $\phi\in U_d$ we set
$$
m_v(\phi)=  \max_{\nu=1,\ldots,n;\kappa=1,\ldots,k}|\phi_{\nu}'(d_{\kappa}(\phi))D_ed_{\kappa}(\phi)|_{L_c(C_n,\mathbb{R})}).	
$$

\begin{proposition} 
Suppose $\tau:V\to C^1$ is continuously differentiable and satisfies $(\tau(y))'(0)=1$ for all $y\in V$, and $R:U\to C^1_n$ is given by
$R_{\nu}(\phi)=g_{\nu}(v(\phi))\tau(v(\phi))$. 

\medskip

(i) Then every derivative $DR(\phi),\,\phi\in U$, has a linear extension $D_eR(\phi):C_n\to C^1_n$ so that the map 
$$
U\times C_n\ni(\phi,\chi)\mapsto D_eR(\phi)\chi\in C^1_n
$$
is continuous.

\medskip

(ii) For all $\phi\in U$ and $y=v(\phi)$, and for all $\chi\in C_n$, we have
\begin{eqnarray*}
|I_nD_eR(\phi)\chi| & \le & kn^2(1+m_v(\phi))(m_{Dg}(y)|I\,\tau(y)|\\
& & +m_g(y)\max_{\mu=1,\ldots,kn}|I\,D_{\mu}\tau(y)1|)\cdot|\chi|.
\end{eqnarray*}
\end{proposition}

{\bf Proof.} 1. On (i). Let $\phi\in U$, $y=v(\phi)$, and $\hat{\phi}\in C^1_n$, $\nu=1,\ldots,n$. Then, 
\begin{eqnarray*}
DR_{\nu}(\phi)\hat{\phi} & = & Dg_{\nu}(y)Dv(\phi)\hat{\phi}\cdot \tau(y)+g_{\nu}(y)D\tau(y)Dv(\phi)\hat{\phi}\\
& = & \sum_{\mu=1}^{kn}\left([Dv_{\mu}(\phi)\hat{\phi}]D_{\mu}g_{\nu}(y)1\cdot\tau(y)+g_{\nu}(y)[Dv_{\mu}(\phi)\hat{\phi}]D_{\mu}\tau(y)1\right).
\end{eqnarray*}
Recall Proposition 2.1. Replacing $Dv_{\mu}(\phi)\hat{\phi}$ in the preceding formula by $D_ev_{\mu}(\phi)\chi$ with $\chi\in C_n$ yields linear extensions
$D_eR_{\nu}(\phi):C_n\to C^1$ for each $\phi\in U$. Obtain $D_eR(\phi):C_n\to C^1_n$ by taking the maps $D_eR_{\nu}(\phi)$ as its components. For $\phi\in U,\,\chi\in C_n,\,\nu=1,\ldots,n$ , we have
\begin{eqnarray*}
I\,D_eR_{\nu}(\phi)\chi & = & \sum_{\mu=1}^{kn}\left([D_ev_{\mu}(\phi)\chi]D_{\mu}g_{\nu}(y)1\cdot\,I\,\tau(y)+\right.\\
& & \left.g_{\nu}(y)[D_ev_{\mu}(\phi)\chi]\,I\,D_{\mu}\tau(y)1\right).
\end{eqnarray*}
From the formula for $D_ev_{\mu}(\phi)\chi$  according to Proposition 2.1 we see that for all $\phi\in U$, $\chi\in C_n$, and $\mu\in\{1,\ldots,kn\}$,
$$
|D_ev_{\mu}(\phi)\chi|\le|\chi|+|\phi_{\nu}'(d_{\kappa}(\phi))D_ed_{\kappa}(\phi)\chi|\le
(1+m_v(\phi))|\chi|,
$$
with $\kappa\in\{1,\ldots,k\}$ and $\nu\in\{1,\ldots,n\}$ given by $\mu=(\kappa-1)n+\nu$. It follows that for all $\phi\in U$ and $\chi\in C_n$, and $\nu\in\{1,\ldots,n\}$,
\begin{eqnarray*}
|I\,D_eR_{\nu}(\phi)\chi| & \le & \sum_{\mu=1}^{kn}\left((1+m_v(\phi))|\chi|m_{Dg}(y)|I\,\tau(y)|+\right. \\
& & \left.m_g(y)(1+m_v(\phi))|\chi||I\,D_{\mu}\tau(y)1|\right),
\end{eqnarray*}
which yields
\begin{eqnarray*}
|I_nD_eR(\phi)\chi| & \le & \sum_{\nu=1}^n|I\,D_eR_{\nu}(\phi)\chi|\\
& \le & kn^2(1+m_v(\phi))(m_{Dg}(y)|I\,\tau(y)|+\\
& & m_g(y)\max_{\mu=1,\ldots,kn}|I\,D_{\mu}\tau(y)1|)|\chi|.\quad\Box
\end{eqnarray*}

\medskip

 We fix some $\epsilon\in(0,1)$ and apply Lemma 2.3 for ${\mathcal V}=V$ and $H=H_c$, $c>0$, given by
$$
H_c(y)=\frac{\epsilon}{kn^2(1+c)(1+m_g(y)+m_{Dg}(y))}.
$$
With the resulting map $\tau=\tau_c$, which satisfies $(\tau_c(y))'(0)=1$ for all $y\in V$,  we define $A_c:U\to C^1_n$ by $A_c(\phi)=\phi-R_c(\phi)$ with 
$$ R_{c,\nu}(\phi)=f_{\nu}(\phi)\tau_c(v(\phi))=g_{\nu}(v(\phi))\tau(v(\phi))
$$
for $\nu=1,\ldots,n$ and $\phi\in U$.
Let
$$
U_c=\{\phi\in U:m_v(\phi)<c\}.
$$

\begin{corollary} 
(Smallness of the remainder)
Let $c>0$ be given. Then
$$
|I_nR_c(\phi)|<\epsilon\,\,\mbox{for all}\,\,\phi\in U
$$
and
$$
|I_nD_eR_c(\phi)|_{L_c(C_n,C_n)}\le\epsilon\,\mbox{for all}\,\,\phi\in U_c.
$$	
\end{corollary}

{\bf Proof.}  Using Proposition 3.1 and Lemma 2.3 we infer
$$
|I_nR_c(\phi)|\le n\,m_g(v(\phi))|I\tau_c(v(\phi))|\le n\,m_g(v(\phi))H_c(v(\phi))<\epsilon
$$
for all $\phi\in U$, and for each $\phi\in U_c$,
\begin{eqnarray*}
|I_nD_eR_c(\phi)\chi| & \le & kn^2(1+c)(m_{Dg}(v(\phi))H_c(v(\phi))+\\
& &
m_g(v(\phi))H_c(v(\phi)))|\chi|\le\epsilon|\chi|\,\,\mbox{for all}\,\,\chi\in C_n.\quad\Box
\end{eqnarray*}

\medskip

Obviously, $U=\cup_{c>0}\,U_c$ and $U_c\supset U_{c'}$ for $c>c'$. As $D_ed_{\kappa}(\phi)$ is in general not continuous with respect to $\phi\in U_d$ the sets $U_c$ do not need to be open. 
But we have
$$
U=\cup_{c>0}\,int\,U_c.
$$
Let us show $U\subset\cup_{c>0}\,int\,U_c$. Let $\phi\in U$ be given.  Continuity of differentiation $C^1\to C$ and of evaluation $C\times[-h,0]\to\mathbb{R}$ and of all $d_{\kappa}$ combined yield  that each map $U\ni\psi\mapsto\psi_{\nu}'(d_{\kappa}(\psi))\in\mathbb{R}$ is continuous, hence locally bounded. Also each map $U\ni\psi\mapsto D_ed_{\kappa}(\psi)\in L_c(C_n,\mathbb{R})$ is locally bounded. It follows that there exist a neighbourhood $N\subset U$ of $\phi$ and $c>0$ with $m_v(\phi)<c$ on $N$. Hence $N\subset U_c$ which means $\phi\in\,int\,U_c$. $\Box$

\begin{corollary} 
Each derivative $DA_c(\phi):C^1_n\to C^1_n$, $c>0$ and $\phi\in U_c$, is a topological isomorphism. The restriction of $A_c$ to $\,int\,U_{c}$ is an open mapping which defines an open mapping from the submanifold $X_f\cap\,int\, U_c$ into the closed subspace $X_0$, with
$A_c(\phi)=\phi$ on $(X_f\cap\,int\, U_c)\cap X_0$.
\end{corollary}

{\bf Proof.} 1. Let $c>0$ and $\phi\in U_c$ be given.

\medskip

1.1 On injectivity of $DA_c(\phi)$. For $\tilde{\phi}\in C^1_n$
with $DA_c(\phi)\tilde{\phi}=0$ we have 
$$
\tilde{\phi}=DR_c(\phi)\tilde{\phi}=D_eR_c(\phi)I_n\tilde{\phi}.
$$
With Corollary 3.2,
$$
|I_n\tilde{\phi}|=|I_nD_eR_c(\phi)I_n\tilde{\phi})|\le\epsilon
|I_n\tilde{\phi}|.
$$
Using $0<\epsilon<1$ we get $I_n\tilde{\phi}=0$, hence $\tilde{\phi}=0$.

\medskip

1.2 On surjectivity of $DA_c(\phi)$. Let $\psi\in C^1_n$ be given, Using Corollary 3.2 and $0<\epsilon<1$ we obtain that the map $C_n\ni\chi\mapsto I_n\psi+I_nD_eR_c(\phi)\chi\in C_n$
is a contraction. Its unique fixed point $\chi_{\ast}$ belongs to $C^1_n$
since $\chi_{\ast}=I_n(\psi+D_eR_c(\phi)\chi_{\ast})$ with $\psi\in C^1_n$ and $D_eR_c(\phi)\chi_{\ast}\in C^1_n$. It follows that
$$
C^1_n\ni\chi_{\ast}=\psi+D_eR_c(\phi)\chi_{\ast}=\psi+DR_c(\phi)\chi_{\ast},
$$
hence 
$$
DA_c(\phi)\chi_{\ast}=\chi_{\ast}-DR_c(\phi)\chi_{\ast}=\psi.
$$

\medskip

1.3 By Parts 1.1 and 1.2, the continuous linear map $DA_c(\phi)$ is an isomorphism $C^1_n\to C^1_n$. Due to the Open Mapping Theorem its inverse also is continuous. So $DA_c(\phi)$ is a topological isomorphism.

\medskip

2. The result of Part 1 allows us to apply the local Inverse Mapping Theorem at each $\phi\in\, int\,U_c$. It follows easily that the restriction of $A_c$ to $int\,U_c$ is an open mapping.

\medskip

3. Recall that $A_c$ maps $X_f\cap\,int\,U_c$ into $X_0$, with $A_c(\phi)=\phi$ on $(X_f\cap\,int\, U_c)\cap X_0$. We abbreviate $A=A_c$ and $X=X_f\cap\,int\,U_c$ and consider the induced map $\alpha:X\ni\phi\mapsto A(\phi)\in X_0$. Proof that the derivatives $D\alpha(\phi)$, $\phi\in X$, are topological isomorphisms: Let $\phi\in X$ be given. Then 
$$
D\alpha(\phi)\psi=DA(\phi)\psi\in T_{A(\phi)}X_0=X_0\,\,\mbox{for all}\,\,\psi\in T_{\phi}X.
$$ 
It follows that $D\alpha(\phi):T_{\phi}X\to T_{\alpha(\phi)}X_0$ is injective and continuous. In order to obtain surjectivity, use that
$DA(\phi):C^1_n\to C^1_n$ is an isomorphism, whereby 
$DA(\phi)T_{\phi}X\subset X_0$ has the same codimension $n$ in $C^1_n$ as $T_{\phi}X$. The previous inclusion and the fact that also $X_0$ has codimension $n$ in $C^1_n$ combined yield $DA(\phi)T_{\phi}X=X_0$, hence
$$
D\alpha(\phi)T_{\phi}X=T_{\alpha(\phi)}X_0.
$$
Due to the Open Mapping Theorem $(D\alpha(\phi))^{-1}$ is continuous.

\medskip  

4. We infer that $A$ defines an open mapping from the submanifold $X$
into the submanifold $X_0$. $\Box$

\begin{theorem} 
For each $c>0$ the map $A_c$ is injective on every convex open subset $\hat{U}_c\neq\emptyset$ of $U_c$ and defines a diffeomorphism $\hat{A}_c$ from $\hat{U}_c$ onto an open subset of $C^1_n$ which maps $X_f\cap\hat{U}_c$ onto an open subset of the subspace $X_0$, with $\hat{A}_c(\phi)=\phi$ on $(X_f\cap\hat{U}_c)\cap X_0$.
\end{theorem}

{\bf Proof.} 1. (Injectivity on $\hat{U}_c$) Let $\phi,\chi$ in $\hat{U}_c$ be given with $A_c(\phi)=A_c(\chi)$. We need to show $\phi=\chi$. By convexity,
$$
\phi+[0,1](\chi-\phi)\subset\hat{U}_c.
$$
Notice that $\chi-\phi=R_c(\chi)-R_c(\phi)$. It follows that
\begin{eqnarray*}
	I_n(\chi-\phi) & = &  I_n(R_c(\chi)-R_c(\phi))
	=I_nR_c(\chi)-I_nR_c(\phi)\\
	& = & \int_0^1D(I_nR_c)(\phi+t(\chi-\phi))(\chi-\phi)dt\\
	& = & \int_0^1I_nDR_c(\phi+t(\chi-\phi))(\chi-\phi)dt\\
	& = & \int_0^1I_nD_eR_c(\phi+t(\chi-\phi))(I_n(\chi-\phi))dt.
\end{eqnarray*}
Using Corollary 3.2 we get
$|I_n(\chi-\phi)|\le\epsilon|I_n(\chi-\phi)|$, hence
$I_n(\chi-\phi)=0$, and thereby $\chi-\phi=0$. 

\medskip

2. Due to Corollary 3.3 all derivatives $DA_c(\phi)$, $\phi\in\hat{U}_c\subset\,int\,U_c$, are topological isomorphisms, which implies that on $\hat{U}_c$ the map $A_c$ is locally invertible by continuously differentiable maps. In combination with injectivity of $A_c$ on $\hat{U}_c$ this yields that the inverse on the open set $A_c(\hat{U}_c)$ is continuously differentiable. For the remaining part of the assertion, use Corollary 3.3. $\Box$ 

\section{Smoothness replacing convexity} 

In this section we consider maps $g:\mathbb{R}^{kn}\supset V\to\mathbb{R}^n$ and
$\Delta:C_n\supset W_{\Delta}\to[-h,0]^k\subset\mathbb{R}^k$, $V$ and $W_{\Delta}$ open, and assume that both maps are continuously differentiable. 

\medskip

Then $U_d=I_n^{-1}(W_{\Delta})=W_{\Delta}\cap C^1_n$  is open in $C^1_n$, and
$d:U_d\to[-h,0]^k\subset\mathbb{R}^k$ given by $d(\phi)=\Delta(I_n\phi)$
is continuously differentiable and has property (e), with $D_ed(\phi)=D\Delta(I_n\phi)\in L_c(C_n,\mathbb{R}^k)$ for $\phi\in U_d$, see Part (i) of the Remarks at the begin of Section 2. So the hypothesis (g,d) from Section 1 is satisfied by $g$ and $d$.

\medskip

The map $w:W_{\Delta}\to\mathbb{R}^{kn}$ given by 
$$
w_{\mu}(\chi)=\chi_{\nu}(\Delta_{\kappa}(\chi))\,\,\mbox{for}\,\,\mu\in\{1,\ldots,kn\},
$$
with $\kappa\in\{1,\ldots,k\}$ and $\nu\in\{1,\ldots,n\}$ determined by $\mu =(\kappa-1)n+\nu$,
is continuous but in general not locally Lipschitz, let alone differentiability.

\medskip

We assume in addition that the open subset $W=W_{\Delta}\cap w^{-1}(V)$ of $C_n$ is non-empty and $W\neq C_n$.

\medskip

For $v=w|U_d$, $v(\phi)=\phi(d(\phi))$, so Proposition 2.1 applies and yields that $v$ is continuously differentiable with property (e). 

\medskip 

As in Section 2 we consider the open set 
$$
U=U_d\cap v^{-1}(V)=(W_{\Delta}\cap C^1_n)\cap (w^{-1}(V)\cap C^1_n)=W\cap C^1_n\subset C^1_n.
$$

\medskip

Proposition 2.2 applies and yields that $f:C^1_n\supset U\to\mathbb{R}^n$ given by $f(\phi)=g(v(\phi))$ is continuously differentiable with property (e), and that $X_f$ if non-empty is a continuously differentiable submanifold of codimension $n$ in $C^1_n$.

\medskip

We turn to subsets which exhaust $W$ and $U$, respectively, and begin with $W$: For $q>0$ let
$$
W_q=\{\chi\in W:|D\Delta_{\kappa}(\chi)|_{L_c(C_n,\mathbb{R})}<q\,\,\mbox{for}\,\,\kappa=1,\ldots,k\}.
$$
The sets $W_q$ are open and $W=\cup_{q>0}W_q$. The preimages
\begin{eqnarray*}
U_q & = & I_n^{-1}(W_q)\\
& = & \{\phi\in I_n^{-1}(W):|D\Delta_{\kappa}(I_n\phi)|_{L_c(C_n,\mathbb{R})}<q\,\,\mbox{for}\,\,\kappa=1,\ldots,k\}\\
& = & \{\phi\in I_n^{-1}(W):|D_ed_{\kappa}(\phi)|_{L_c(C_n,\mathbb{R})}<q\,\,\mbox{for}\,\,\kappa=1,\ldots,k\}
\end{eqnarray*}
for $q>0$ are open, and $U=\cup_{q>0}U_q$. For $q>0$ and  $b>0$ we consider the open sets
$$
U_{q,b}=\{\phi\in U_q:|\partial\phi|<b\},
$$
which satisfy $U=\cup_{q>0,b>0}U_{q,b}$, and
$$
U_{q,b}\subset U_c
$$
for $U_c$ from Section 3 with $c=bq$.

For $q>0,b>0,$ and $0<\delta<1$, we define the map 
$H_{q,b,\delta}:V\to(0,\infty)\subset\mathbb{R}$ by
$$
H_{q,b,\delta}(y)=\frac{\delta}{kn^2(1+bq)(1+m_g(y)+m_{Dg}(y))}
$$
and apply Lemma 2.3 with ${\mathcal V}=V$ and $H=H_{q,b,\delta}$. With the resulting map $\tau=\tau_{q,b,\delta}$ we define
$$
A_{q,b,\delta}:U_{q,b}\ni\phi\mapsto \phi-R_{q,b,\delta}(\phi)\in C^1_n
$$
by
$$ R_{q,b,\delta,\nu}(\phi)=f_{\nu}(\phi)\tau_{q,b,\delta}(v(\phi))=g_{\nu}(v(\phi))\tau_{q,b,\delta}(v(\phi))
$$
for $\nu=1,\ldots,n$ and $\phi\in U_{q,b}$. The proofs of Proposition 3.1,  Corollary 3.2, and Corollary 3.3, remain valid, with $R_{q,b,\delta}$ in place of $R$, and $U_{q,b}$ in place of each of $U$, $U_c$, and $int\,U_c$, and $\delta$ in place of $\epsilon$. In particular, for all $q>0$, $b>0$, and $\phi\in U_{q,b}$,
$$
|I_nR_{q,b,\delta}(\phi)|<\delta\quad\mbox{and}\quad|I_nD_eR(\phi)|_{L_c(C_n,C_n)}\le\delta.
$$
Due to Corollary 3.3 each derivative $DA_{q,b,\delta}(\phi)$, $\phi\in U_{q,b}$, is a topological isomorphism and $A_{q,b,\delta}$ is an open map which induces an open map from the submanifold $X_f\cap U_{q,b}$ into the subspace $X_0$.

\begin{proposition} 
For all $q>0,b>0,$ and $\delta\in(0,1)$ the restriction $A_{q,b,\delta,\ast}$ of $A_{q,b,\delta}$ to 	
$$
U_{q,b,\delta}=\{\phi\in U_{q,b}:dist(I_n\phi,C_n\setminus W_q)>2\delta\}	
$$ 
is injective.
\end{proposition}

The sets $U_{q,b,\delta}$ are open and $U=\cup_{q>0,b>0,\delta>0}U_{q,b,\delta}$. For $\phi\in U_{q,b,\delta}$, $N_{2\delta}(I_n\phi)\subset W_q$.

\medskip

{\bf Proof} of Proposition 4.1. We abbreviate $A=A_{q,b,\delta}$ and $R=R_{q,b,\delta}$. Assume $A(\phi)=A(\psi)$ for 
$\phi$ and $\psi$ in $U_{q,b,\delta}$. Let $\zeta=A(\phi)$.Then
$$
|I_n\phi-I_n\zeta|=|I_n\phi-I_nA(\phi)|=|I_nR(\phi)|<\delta.
$$
It follows that 
$$
I_n\phi\in N_{\delta}(I_n\zeta)\subset N_{2\delta}(I_n\phi)\subset W_q.
$$
Analogously, $I_n\psi\subset N_{\delta}(I_n\zeta)\subset W_q$. 
By convexity, $I_n\phi+[0,1](I_n\psi-I_n\phi)\subset N_{\delta}(I_n\zeta)\subset W_q$. It follows that
$$
\phi+[0,1](\psi-\phi)\subset I_n^{-1}(N_{\delta}(I_n\zeta))\subset I_n^{-1}(W_q)=U_q\subset C^1_n.
$$
From $|\partial\phi|<b$ and $|\partial\psi|<b$,
$|\partial(\phi+t(\chi-\phi))|=|t\partial\chi+(1-t)\partial\phi|<b$ for all $t\in[0,1]$, and we deduce
$$
\phi+[0,1](\chi-\phi)\subset U_{q,b}.
$$
The proof of $\phi=\psi$ is completed as in Part 1 of the proof of 
Theorem 3.4, with $\hat{U}_c$ replaced by $U_{q,b}$ and $R_c$ replaced by $R=R_{q,b,\delta}$. $\Box$

\medskip

Arguing as in the proof of Theorem 3.4 we obtain the following result.

\begin{theorem} 
Assume that $g:\mathbb{R}^{kn}\supset V\to\mathbb{R}^n$ and
$\Delta:C_n\supset W_{\Delta}\to[-h,0]^k\subset\mathbb{R}^k$, $V$ and $W_{\Delta}$ open, are continuously differentiable, and $\emptyset\neq W\neq C_n$ for $W=W_{\Delta}\cap w^{-1}(V)$ with $w:W_{\Delta}\to\mathbb{R}^{kn}$ given by 
$$
w_{\mu}(\chi)=\chi_{\nu}(d_{\kappa}(\chi))\,\,\mbox{for}
\,\,\mu\in\{1,\ldots,kn\}
$$
where $\kappa\in\{1,\ldots,k\}$ and $\nu\in\{1,\ldots,n\}$ determined by $\mu =(\kappa-1)n+\nu$.
Then for all $q>0,b>0,$ and $\delta\in(0,1)$ with $U_{q,b,\delta}\neq\emptyset$,
there exists a diffeomorphism $A_{q,b,\delta,\ast}:U_{q,b,\delta}\to C^1_n$
onto an open subset of $C^1_n$ which maps $X_f\cap U_{q,b,\delta}$ onto an open subset of the space $X_0$ and satisfies $A_{q,b,\delta,\ast}(\phi)=\phi$ on $(X_f\cap U_{q,b,\delta})\cap X_0$.
\end{theorem}

\section{Examples}

The main results from \cite{KR,W6,KW2} on graph and almost graph representations of solution manifolds imply that as in Theorems 3.4 and  4.2  diffeomorphisms transform the solution manifolds to open subsets of the space $X_0$, under different hypotheses on delay functionals and feedback maps in the system (2), or more generally on the functional $f$ in Eq. (1). In this section we give examples of $g:\mathbb{R}\supset V\to\mathbb{R}$ and $d:C^1\supset U_d\to[-h,0]\subset\mathbb{R}$ so that Theorems 3.4 or 4.2 apply to the solution manifold $X_f\subset C^1$
associated with Eq. (1) for $f:C^1\supset U\to\mathbb{R}$ given by 
$$
U=\{\phi\in U_d:\phi(d(\phi))\in V\}\quad\mbox{and}\quad f(\phi)=g(\phi(d(\phi))),
$$
and yield regularization for open subsets $X_f\cap\hat{U}_c$ (Theorem 3.4) or $X_f\cap U_{q,b,\delta}$ (Theorem 4.2) whereas for these subsets all of the corresponding theorems from \cite{KR,W6,KW2} fail since at least one of their hypotheses is violated. These hypotheses are the following:  

\medskip

(I) {\it The delay functional $d$ factorizes
\begin{equation*}
	d(\phi)=\delta(LI\phi),\quad\phi\in U_d,
\end{equation*}
into a continuous linear map $L:C\to F$ and a continuously differentiable function $\delta:F\supset W\to[-h,0]\subset\mathbb{R}$, where $F$ is a finite-dimensional normed vector space, $W\subset F$ open, and $LIU_d\subset W$.}

\medskip

(Required in \cite[Theorem 5.1]{W6} and \cite[Theorem 3.5]{KW2}.)

\medskip

(II) {\it Delays are bounded away from zero in the sense that there exists $s\in(-h,0)$
so that $f(\phi)=f(\psi)$ for all $\phi,\psi$ in $U$
with $\phi(t)=\psi(t)$ for all $t\in[-h,s]\cup\{0\}$.}

\medskip

(Required in \cite[Theorem 2.4]{W6}.)

\medskip

(III) {\it The set of extended derivatives $D_ef(\phi)$, $\phi\in U$, is bounded in $L_c(C,\mathbb{R})$.}

\medskip

(Required in the proof of \cite[Lemma 1]{KR}.)

\medskip

Notice that the factorization property (I) implies that $d$ is constant on the intersection of $U_d$ with the vector space $(LI_n)^{-1}(0)\subset C^1$.

\medskip

In order to specify $d$ we choose continuously  differentiable functions $p:\mathbb{R}\to[0,\infty)\subset\mathbb{R}$ and $\eta:\mathbb{R}\to[-h,0]\subset\mathbb{R}$ with bounded derivatives and
$$
p(\xi)=0\,\,\mbox{on}\,\,(-\infty,0],\,\,p'(\xi)>0\,\,\mbox{on}\,\,(0,\infty),\,\,0<\eta'(\xi)\,\,\mbox{everywhere}.
$$
We define $\Delta:W_{\Delta}\to[-h,0]\subset\mathbb{R}$ and $d:U_d\to[-h,0]\subset\mathbb{R}$ by
$$
W_{\Delta}=C\quad\mbox{and}\quad\Delta(\chi)=\eta\left(\int_{-h}^0p\circ\chi\right),
$$
$$
U_d=I^{-1}(W_{\Delta})=C^1,\quad d(\phi)=\Delta(I\phi).
$$

\begin{corollary} 
(i) $\Delta$ and $d$ are continuously differentiable and $d$ has property (e).

\medskip

(ii) There exists $q>0$ with $|D\Delta(\chi)|_{L_c(C,\mathbb{R})}\le q$
for all $\chi\in W_{\Delta}=C$ and $|D_ed(\phi)|_{L_c(C,\mathbb{R})}\,\le\,q$ for all $\phi\in U_d=C^1$.

\medskip

(iii) For every open set ${\mathcal O}\subset U_d=C^1$ with $0\in{\mathcal O}$ and for every subspace 
$Z\subset C^1$ with $Z\neq\{0\}$ the delay functional $d$ is not constant on ${\mathcal O}\cap Z$.
\end{corollary}

{\bf Proof} $\,\, 1.$ The continuous differentiability of $\Delta$ follows by means of the chain rule from the continuous differentiability of the substitution operator $C\ni\chi\mapsto p\circ\chi\in C$ (compare \cite[Lemma 1.5 in Appendix IV]{DvGVLW}) in combination with continuity of the linear map $C\ni\chi\mapsto\int_{-h}^0\chi\in\mathbb{R}$ and continuous differentiability of $\eta$. We obtain
$$
D\Delta(\chi)\hat{\chi}=\eta'\left(\int_{-h}^0p\circ\chi\right))\int_{-h}^0p'(\chi(t))\hat{\chi}(t)dt
$$
for all $\chi,\hat{\chi}$ in $C$. By the chain rule also $d=\Delta\circ I$ is continuously differentiable, and due to assertion (i) of the Remark in Section 2 $d$ has property (e) with $D_ed(\phi)=D\Delta(I\phi)$ on $U_d=C^1$.

\medskip

2. Assertion (ii) follows from the boundedness of $\eta'$ and $p'$  by means of the formulae for
$D\Delta$ and $D_ed$ in Part 1.

\medskip

3. On (iii). Choose $\phi\in ({\mathcal O}\cap Z)\setminus\{0\}$. As ${\mathcal O}$ is a neighbourhood of $0\in C^1$ multiplication by a suitably small number $\epsilon\neq0$ yields $\epsilon\phi\in {\mathcal O}\cap Z$, $2 \epsilon\phi\in {\mathcal O}\cap Z$ and $\epsilon\phi(s)>0$ for some $s\in(-h,0)$. We infer
\begin{eqnarray*}
\int_{-h}^0p\circ(\epsilon\phi) & = & 
\int_{\{t\in[-h,0]:\phi(t)>0\}}p\circ(\epsilon\phi)\\
& < & \int_{\{t\in[-h,0]:\phi(t)>0\}}p(2\epsilon\phi(t))dt=\int_{-h}^0p(2\epsilon\phi(t))dt,
\end{eqnarray*}
and injectivity of $\eta$ yields
$$
d(\epsilon\phi)=\eta\left(\int_{-h}^0p\circ(\epsilon\phi)\right)\neq
\eta\left(\int_{-h}^0p(2\epsilon\phi(t))dt\right)=d(2\epsilon\phi).\quad \Box
$$

Let $g:V\to\mathbb{R}$ be a continuously differentiable function on $V=(-\infty,\gamma)$ with $\gamma>0$ which is injective and satisfies $|g'(-n)|\to\infty$ as $n\to\infty$. Consider $w:U_{\Delta}\to\mathbb{R}$ and $v:U_d\to\mathbb{R}$ as chosen before Corollary 5.1. Then $W=W_{\Delta}\cap w^{-1}(V)$ equals $\{\chi\in C:\chi(\Delta(\chi))<\gamma\}$ and satisfies $\emptyset\neq W\neq C$ as it is required in Section 4, and
$U=U_d\cap v^{-1}(V)$ equals $\{\phi\in C^1:\phi(d(\phi))<\gamma\}$. 
Recall that as in Section 4 $f:U\to\mathbb{R}$ given by $f(\phi)=g(v(\phi))$ is continuously differentiable with property (e).

\medskip

For $n\in\mathbb{N}$ let ${\bf n}$ denote the constant function $[-h,0]\to\mathbb{R}$ with value $n$.

\begin{corollary} 
(i) For every neighbourhood $N\subset U$ of $0$ in $C^1$ and for every $s\in(-h,0)$ there exist $\phi$ and $\psi$ in $N$ such that
$\phi(t)=\psi(t)$ on $[-h,s]\cup\{0\}{\tiny }$ and $f(\phi)\neq f(\psi)$,

\medskip

(ii) $|D_ef(-{\bf n})|_{L_c(C,\mathbb{R})}|\,\,\to\,\,\infty\,\,\mbox{as}\,\,n\to\infty.$
\end{corollary}

{\bf Proof} $\,\,$ 1. On (i). Obviously, $0\in U$. Let a neighbourhood $N\subset U$ of $0$ in $C^1$ and $s\in(-h,0)$ be given. Choose strictly increasing $\phi,\psi$ in $N$ with $\phi(t)=\psi(t)$ on $[-h,s]\cup\{0\}$ and $\phi(s)=0=\psi(s)$, and $\phi(t)<\psi(t)$ on $(s,0)$. By the properties of $p$,
$$
\int_{-h}^0p\circ\phi=\int_s^0p\circ\phi<\int_s^0p\circ\psi=\int_{-h}^0p\circ\psi.
$$
Using $\eta'(\xi)>0$ everywhere we get 
$$
d(\phi)=\eta\left(\int_{-h}^0p\circ\phi\right)<\eta\left(\int_{-h}^0p\circ\psi\right)=d(\psi).
$$ 
As $\phi$ is strictly increasing, and $\phi\le\psi$, we obtain
$$
\phi(d(\phi))<\phi(d(\psi))\le\psi(d(\psi)),
$$
and the injectivity of $g$ yields
$$
f(\phi)=g(\phi(d(\phi)))\neq g(\psi(d(\psi)))=f(\psi)).
$$
2. On (ii). For each $n\in\mathbb{N}$ we have $-{\bf n}\in U$ and 
\begin{eqnarray*}
|D_ef(-{\bf n})|_{L_c(C,\mathbb{R})} & \ge & |D_ef(-{\bf n}){\bf 1}|=|Df(-{\bf n}){\bf 1}|=|Dg(v(-{\bf n}))Dv(-{\bf n}){\bf 1}|\\
& = & 
|g'(-n)({\bf 1}(d(-{\bf n}))+(-{\bf n})'(d(-{\bf n}))Dd(-{\bf n}){\bf 1})|\\
& = & |g'(-n)(1+0)|=|g'(-n)|,
\end{eqnarray*}
hence $|D_ef(-{\bf n})|_{L_c(C,\mathbb{R})}\to\infty$ as $n\to\infty. \quad \Box$

\medskip

In order to verify the hypotheses of Theorem 3.4 we choose $b>0$ and $c>bq$ with $q$ according to Corollary 5.1 (ii). Then $m_v(\psi)<c$ for all $\psi\in U=\{\phi\in C^1:\phi(d(\phi))<\gamma\}$ with $|\partial\psi|<b$, and the open convex neighbourhood
$$
\hat{U}_c=\{\phi\in C^1:\phi(t)<\gamma\,\,\mbox{on}\,\,[-h,0],\,|\partial\phi|<b\}
$$
of $0\in C^1$ is contained in $U_c=\{\phi\in U:m_v(\phi)<c\}$.
Therefore Theorem 3.4 applies to $\hat{U}_c$ and yields regularization of $X_f\cap\hat{U}_c$.

\medskip

Corollary 5.1 (iii) in combination with the remark following condition (III) shows that the restriction of $f$ to the open subset $\hat{U}_c\ni 0$ of $U_d=C^1$ violates condition (I). As $\hat{U}_c\subset U$ is an open neighbourhood of $0$ in $C^1$ Corollary 5.2 (i) shows that condition (II) is violated. As $-{\bf n}\in\hat{U}_c$ for all $n\in\mathbb{N}$ also condition (III) is  violated,

\medskip

Application of Theorem 4.2: Let $w:W_{\Delta}\to\mathbb{R}$ be given by $w(\chi)=\chi(\Delta(\chi))$. For $W=W_{\Delta}\cap w^{-1}(V)=\{\chi\in C:\chi(\Delta(\chi))<\gamma\}$ the hypothesis $\emptyset\neq W\neq C$ from Section 4 is satisfied. Due to Corollary 5.1 (ii) there exists $q>0$ with $|D\Delta(\chi)|_{L_c(C,\mathbb{R})}<q$ for all $\chi\in W_{\Delta}$. Hence the set $W_q=\{\chi\in W:|D\Delta(\chi)|_{L_c(C,\mathbb{R})}<q\}$ equals $W$, and we get $U=I^{-1}(W)=I^{-1}(W_q)=U_q$. As $0\in W$ and $W$ is open there exists $\delta\in(0,1/2)$ with $2\delta<dist(0,C\setminus W)=dist(0,C\setminus W_q)$. We choose some $b>0$ and consider the open sets
$$ 
U_{q,b}=\{\phi\in U_q:|\partial\phi|<b\}=\{\phi\in U:|\partial\phi|<b\}
$$
and
\begin{eqnarray*}
U_{q,b,\delta} & = & \{\phi\in U_{q,b}:dist(I\phi,C\setminus W_q)>2\delta\}\\
& = &  \{\phi\in U:|\partial\phi|<b,\,\,
dist(I\phi,C\setminus W)>2\delta\}\\
& = & \{\phi\in C^1:\phi(d(\phi))<\gamma,\,\,|\partial\phi|<b,\,\,
dist(I\phi,C\setminus W)>2\delta\}.
\end{eqnarray*}
Obviously, $0\in U_{q,b,\delta}$, so $U_{q,b,\delta}$ is non-empty.  Theorem 4.2 applies to $U_{q,b,\delta}$ and yields regularization of $X_f\cap U_{q,b,\delta}$.  

\medskip

Corollary 5.1 (iii) in combination with the remark following condition (III) shows that the restriction of $f$ to the open subset $U_{q,b,\delta}\ni0$ of $U_d=C^1$ violates condition (I). As $U_{q,b,\delta}\subset U$ is an open neighbourhood of $0$ in $C^1$ Corollary 5.2 (i) shows that condition (II) is violated for $f|U_{q,b,\delta}$. In order to see that condition (III) is violated for $f|U_{q,b,\delta}$ notice that for every $n\in\mathbb{N}$ and for each
$$
\chi\in C\setminus W=\{\chi\in C:\chi(\Delta(\chi))\ge\gamma\}
$$ 
we have
$|-{\bf n}-\chi|\ge1$ and thereby  
$dist(-{\bf n},C\setminus W)\ge1>2\delta$. It follows that
$-{\bf n}\in U_{q,b,\delta}$ for all $n\in\mathbb{N}$, and Corollary 5.2 (ii) shows that condition (III) is violated for $f|U_{q,b,\delta}$.


\end{document}